\pdfoutput=1
\RequirePackage{ifpdf}
\ifpdf 
\documentclass[pdftex]{sigma}
\else
\documentclass{sigma}
\fi

\newcommand{\C}{\mathbb{C}}
\newcommand{\R}{\mathbb{R}}
\newcommand{\Z}{\mathbb{Z}}
\newcommand{\N}{\mathbb{N}}

\numberwithin{equation}{section}

\newtheorem{Theorem}{Theorem}[section]
\newtheorem*{Theorem*}{Theorem}

\newtheorem{Lemma}[Theorem]{Lemma}

\newtheorem{Conjecture}[Theorem]{Conjecture}

\theoremstyle{definition}

\newtheorem{Remark}[Theorem]{Remark}

\begin{document}

\allowdisplaybreaks

\renewcommand{\thefootnote}{}

\renewcommand{\PaperNumber}{120}

\FirstPageHeading

\ShortArticleName{Local 3D Weinstein Conjecture}

\ArticleName{Local 3D Weinstein Conjecture\footnote{This paper is a~contribution to the Special Issue on Geometry and Dynamics in memory of Will Merry. The~full collection is available at \href{https://sigma-journal.com/Merry.html}{https://sigma-journal.com/Merry.html}}}

\Author{Umberto Leone HRYNIEWICZ~$^{\rm a}$ and Patrice LE CALVEZ~$^{\rm b}$}

\AuthorNameForHeading{U.L.~Hryniewicz and P.~Le Calvez}

\Address{$^{\rm a)}$~Lehrstuhl f\"ur Geometrie und Analysis, Fachgruppe Mathematik, RWTH Aachen,\\
\hphantom{$^{\rm a)}$}~Pontdriesch 10-12, 52064 Aachen, Germany}
\EmailD{\mail{hryniewicz@mathga.rwth-aachen.de}}
\URLaddressD{\url{https://www.mathga.rwth-aachen.de/home/home/}}

\Address{$^{\rm b)}$~Sorbonne Universit\'e, Universit\'e Paris Cit\'e, CNRS, IMJ-PRG, 75005 Paris, France}
\EmailD{\mail{patrice.le-calvez@imj-prg.fr}}

\ArticleDates{Received June 18, 2026, in final form November 20, 2026; Published online November 29, 2026}

\Abstract{In this note, we introduce a local version of the 3D Weinstein conjecture, asking whether every compact non-empty isolated invariant set of a Reeb flow on a 3-manifold contains a periodic orbit. Our main result establishes this in the case where the compact isolated invariant set is an embedded torus.}


\Classification{37J06; 37J12; 53E50}

\renewcommand{\thefootnote}{\arabic{footnote}}
\setcounter{footnote}{0}

\section{Introduction}

The goal of this paper is to study invariant sets for Reeb flows in dimension three. We start with some basic notions. Consider a space $M$ with a flow $(t,x) \in\R\times M \mapsto \phi^t(x)\in M$. Let $A \subset M$ be arbitrary. The maximal invariant subset of $A$ is defined as
\begin{equation*}
S(A) = \bigcap_{t\in\R} \phi^t(A).
\end{equation*}
The set $A$ is called invariant if $A=S(A)$. When $M$ is a topological space and the flow is continuous, an invariant set $A \subset M$ is called isolated if $A = S(N)$ for some neighborhood $N$ of~$A$. The notion of isolated invariant sets plays an important role in dynamical systems. Let us state and discuss the main problem we would like to investigate.

\begin{Conjecture}[local 3D Weinstein conjecture]
\label{3DW_conj}
A compact non-empty invariant isolated set of a smooth Reeb flow on a $3$-manifold contains a periodic orbit.
\end{Conjecture}

The 3D Weinstein conjecture is a special case: If the ambient $3$-manifold is closed, then it is a compact isolated invariant set. Some of the fundamental tools to attack the 3D Weinstein conjecture were laid out by Hofer~\cite{93}, who confirmed it for the case $M=S^3$, for the case where~$\pi_2(M)$ is non-trivial, and for the case where the contact structure induced by the contact form defining the Reeb flow is overtwisted. Results by Taubes~\cite{Taubes} based on Seiberg--Witten theory settled the general case.

The main theorem of this note confirms Conjecture~\ref{3DW_conj} for invariant tori.

\begin{Theorem}\label{thm_main_1}
Consider a smooth Reeb flow on a $3$-manifold. Any isolated invariant smoothly embedded torus contains a periodic orbit.
\end{Theorem}

Its contrapositive can be seen as a KAM-type result in disguise, we explain. A smooth flow on a two-torus without periodic orbits or rest points is necessarily a suspension of a circle diffeomorphism. This is well known, but we included a detailed argument for the sake of completeness; see Lemma~\ref{lemma_global_circle_section}. The corresponding circle diffeomorphism must have irrational rotation number. This justifies calling such tori \emph{irrational}. Theorem~\ref{thm_main_1} can be reformulated as follows: \emph{Irrational invariant tori for Reeb flows in $3$D are never isolated}. In KAM theory one typically tries to prove that an irrational invariant torus is not isolated by proving that there are invariant tori nearby. Of course, such a strong conclusion requires stronger assumptions: irrationality is not enough. Theorem~\ref{thm_main_1} does not go as far as providing invariant tori nearby, but this is compensated by making significantly less assumptions: only irrationality is required, i.e., there is no need to look at higher order jets of the flow or to worry about small divisors.\looseness=1

The proof of Theorem~\ref{thm_main_1} presented here is a sequence of reductions to the following special case of the main result of~\cite{LeCY}.

\begin{Theorem}[Le Calvez and Yoccoz~\cite{LeCY}]\label{thm_LeCY}
Consider a germ $f$ of an area- and orientation-preserving surface homeomorphism defined around a fixed point~$p$. If $p$ is isolated as an invariant set, then ${\rm ind}(f^j,p) \leq 0$ for infinitely many iterates $j$.
\end{Theorem}

Here ${\rm ind}(f^j,p)$ denotes the Lefschetz index of $f^j$ at $p$, which is well defined for all $j$ by the isolation hypothesis. In fact, the result in~\cite{LeCY} is more precise: it includes a periodicity statement for the iterates~$j$ such that ${\rm ind}(f^j,p) \leq 0$ and a description of ${\rm ind}(f^j,p)$ for those iterates. This extra information is not needed in our proof of Theorem~\ref{thm_main_1}; the statement presented above as~Theorem~\ref{thm_LeCY} was later also proved by Franks using Conley index theory, see~\cite[Proposition~3.1]{Franks}.

\begin{Remark}
It is, of course, natural to ask whether the higher-dimensional version of Conjecture~\ref{3DW_conj} holds, namely, whether non-empty compact invariant isolated sets of smooth Reeb flows always contain a periodic orbit. At the moment, it is not clear that this is a reasonable conjecture.
\end{Remark}

We close the introduction stating the corresponding conjecture for surface maps.

\begin{Conjecture}\label{maps_conj}
A compact non-empty invariant isolated set of a smooth Hamiltonian diffeomorphism of a symplectic surface contains a periodic orbit.
\end{Conjecture}

\section{Proof of Theorem~\ref{thm_main_1}}

Let $\lambda$ be a contact form on some $3$-manifold. The associated Reeb vector field~$X$ is implicitly defined by equations $i_X{\rm d}\lambda = 0$, $i_X\lambda=1$. Its flow will be denoted by $\phi^t$ and called the Reeb flow. Throughout we assume completeness of the Reeb flow, without loss of generality.

Our argument is indirect. Assume, by contradiction, the existence of a smoothly embedded $2$-torus $\Sigma \subset M$ with the following properties:
\begin{itemize}\itemsep=0pt
\item[(P1)] $\Sigma$ is $\phi^t$-invariant.
\item[(P2)] $\Sigma$ is isolated with respect to $\phi^t$.
\item[(P3)] $\Sigma$ contains no periodic orbits of $\phi^t$.
\end{itemize}

\begin{Lemma}\label{lemma_global_circle_section}
There exists a global section for $\phi^t|_\Sigma$ diffeomorphic to a circle.
\end{Lemma}

This lemma is well known; the argument is included for sake of completeness.

\begin{proof}
The first step is to prove that there exists an embedded circle $C$ transverse to $X$. To see this, consider a vector field $Y$ on $\Sigma$ that is pointwise everywhere linearly independent from $X$. Let $\varphi^t$ denote its flow. Choose a recurrent point $p$ for $\varphi^t$. If $t\mapsto \varphi^t(p)$ is periodic, then we can set $C = \varphi^\R(p)$. Suppose $\varphi^t(p)$ is not a periodic trajectory. Choose, with $r_0>0$ small enough, flow box local coordinates $(x,y) \in (-r_0,r_0)^2$ for $\varphi^t$ centered at $p$, meaning that $p \simeq (0,0)$ and $Y\simeq \partial_y$. Write $X = u \, \partial_x+v\, \partial_y$. Then $u(x,y)\neq0$ for every $(x,y) \in (-r_0,r_0)^2$. Choose $a$ such that $a/2 > |v(0,0)/u(0,0)|$. For each $r>0$, consider the rectangle $Q_r = [-r,r]\times [-ar,ar]$, and for every $x_0 \in [-r,r]$ consider a smooth function $h\colon [-ar,0]\to \R$ satisfying
\begin{itemize}\itemsep=0pt
\item $h(-ar)=x_0$, $h(0)=0$, $\|h\|_{L^\infty} \leq |x_0|$,
\item $h^{(j)}(-ar)=h^{(j)}(0)=0$ for all $j\geq1$, $\|h'\|_{L^\infty} \leq 2/a$.
\end{itemize}
We claim that if $r$ is small enough, then $Q_r \subset (-r_0,r_0)^2$ and for each $x_0 \in [-r,r]$ the path~$\gamma$ from $(x_0,-ar)$ to $(0,0)$ given by $\gamma(s) = (h(s),s)$, $s\in[-ar,0]$, is transverse to $X$. In fact, in local coordinates $X = (u(0,0),v(0,0)) + O(r)$ over $Q_r$. Hence, if $r$ is small enough, then $\|v/u\|_{L^\infty(Q_r)} < a/2$ holds. In particular, if $\dot\gamma(s_0)$ is parallel to $X(\gamma(s_0))$ for some $s_0 \in [-ar,0]$, then $h'(s_0)\neq0$, and we get $|v(\gamma(s_0))/u(\gamma(s_0))| = |1/h'(s_0)| \geq a/2$, a contradiction. Our claim is proved. From now, we fix $r>0$ with these properties.

Since $p$ is a recurrent point, there exists a first time $T>0$ such that $\varphi^t(p)$ hits the transversal arc corresponding to $[-r,r]\times\{-ar\} \subset Q_r$. Let $(x_0,-ar)$ be the first hitting point. Note that~$T$ is also the first positive time such that $\varphi^t(p)$ hits the set corresponding to $[-r,r] \times [-ar,0]$. Consider the arc $\gamma$ from $\varphi^T(p)$ to $p$ given in local coordinates by $s\mapsto(h(s),s)$, $s\in[-ar,0]$, where $h$ has the properties listed above. Then $C := \varphi^{[0,T]}(p) \cup \gamma$ is a smooth embedded circle transverse to $X$. Note that $C$ is essential: if not then, since $X$ is transverse to $C$, $X$ would have a zero in the domain bounded by $C$, contradicting the fact that $X$ is nowhere vanishing along~$\Sigma$.

\centerline{\includegraphics[width=0.6\textwidth]{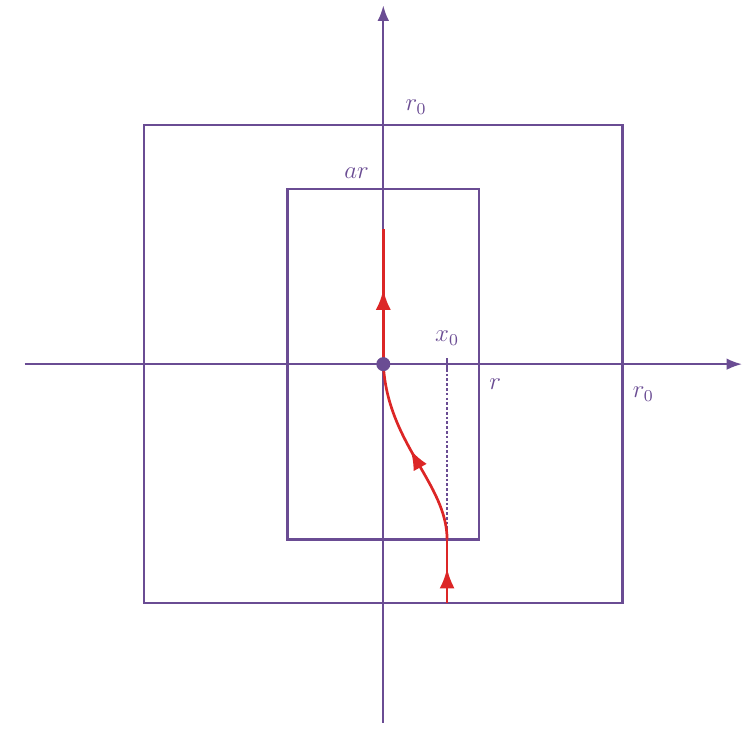}}

The second step is to show that $C$ is a global surface of section for $\phi^t$. To see this we argue by contradiction. Since $C$ is embedded, one can identify $\Sigma$ with $C\times S^1$, where we see the second factor $S^1$ as $[0,1]/\{0\sim1\}$. Then the complement of $C$ gets identified with $C \times (0,1)$. Some trajectory of $\phi^t$ starting at $C \simeq C \times \{0\}$, which will of course immediately enter $C\times(0,1)$, never hits~$C$ again. By transversality of $\phi^t$ to $C$, the $\omega$-limit set of this trajectory is a compact subset of $C \times (0,1)$. The Poincar\'e--Bendixson theorem gives a periodic orbit of~$\phi^t$, in contradiction to~(P3).\looseness=-1
\end{proof}

In view of Lemma~\ref{lemma_global_circle_section}, we can fix a smooth global section $C$ for $\phi^t|_\Sigma$ diffeomorphic to a circle. We can then thicken $C$ in the normal direction of $\Sigma$ to obtain an embedded surface $S$ satisfying $S \pitchfork \Sigma$ and $C = S \cap \Sigma$. Up to shrinking $S$, it can be assumed that $X$ is transverse to $S$. This implies that ${\rm d}\lambda$ induces a smooth area-form $\omega$ on~$S$, and that
\begin{equation*}
R := \bigl\{ p \in S \mid \phi^{(0,+\infty)}(p) \cap S \neq \varnothing \bigr\}
\end{equation*}
is open in $S$.

\centerline{\includegraphics[width=0.45\textwidth]{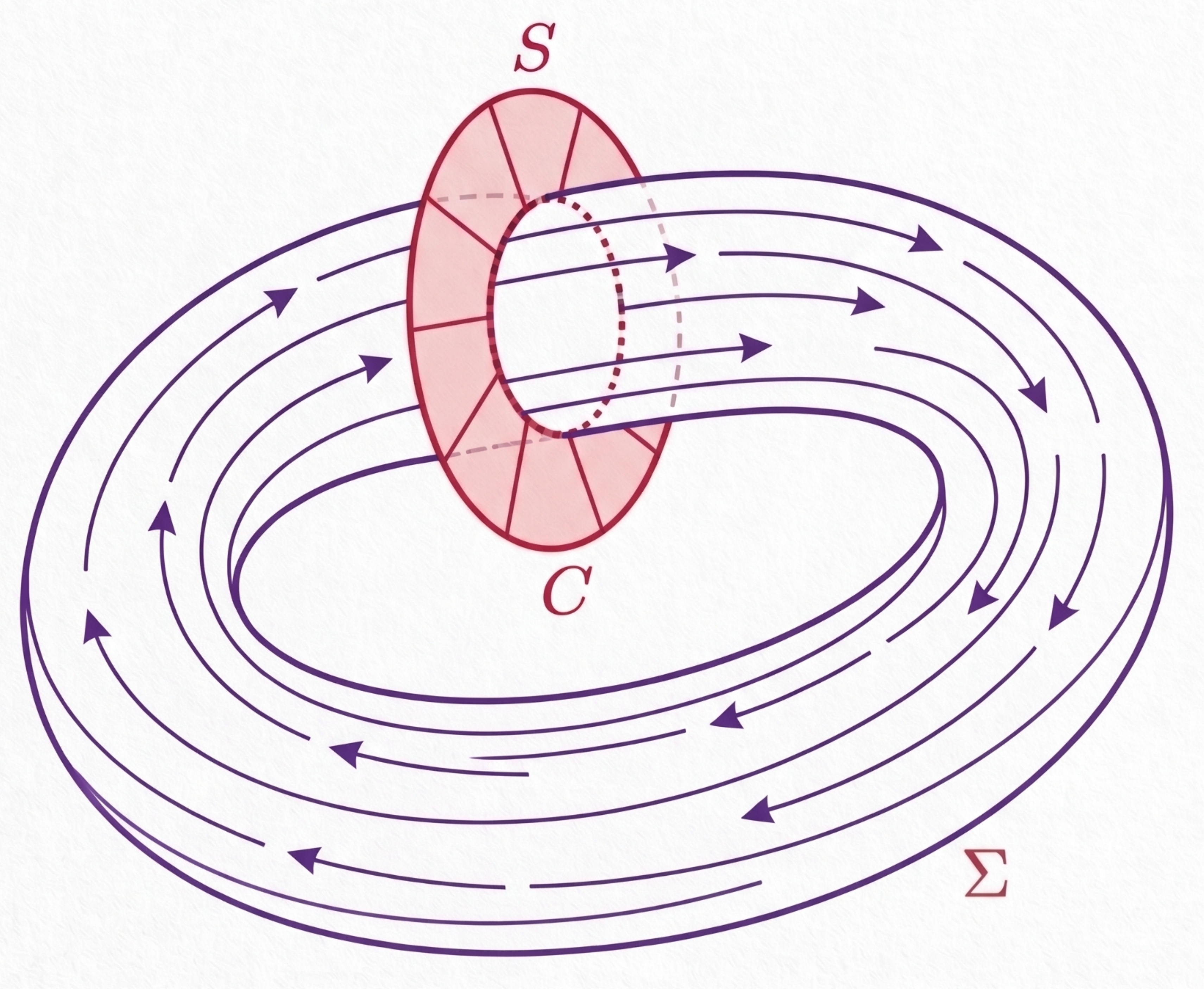}}

Moreover, the first return time function
\begin{equation}\label{return_time_function}
\tau\colon \ R \to (0,+\infty), \qquad \tau(p) = \inf \bigl\{ t > 0 \mid \phi^t(p) \in S \bigr\}
\end{equation}
is smooth. The \emph{first return map}
\begin{equation}\label{return_map}
\psi\colon \ R \to S, \qquad \psi(p) = \phi^{\tau(p)}(p)
\end{equation}
is an $\omega$-preserving embedding onto an open subset of $S$. In fact, choose $p\in R$ and $v_1,v_2\in T_pS$ arbitrarily, and compute
\[ \begin{aligned}
(\psi^*\omega)_p(v_1,v_2)
&= {\rm d}\lambda_{\psi(p)}({\rm d}\psi_p(v_1),{\rm d}\psi_p(v_2)) \\
&= {\rm d}\lambda_{\psi(p)}\bigl({\rm d}\tau_p(v_1)X_{\psi(p)} + {\rm d}\bigl(\phi^{\tau(p)}\bigr)_p(v_1),{\rm d}\tau_p(v_2)X_{\psi(p)} + {\rm d}\bigl(\phi^{\tau(p)}\bigr)_p(v_2)\bigr) \\
&= {\rm d}\lambda_{\phi^{\tau(p)}(p)}\bigl({\rm d}\bigl(\phi^{\tau(p)}\bigr)_p(v_1),{\rm d}\bigl(\phi^{\tau(p)}\bigr)_p(v_2)\bigr)
= \bigl(\bigl(\phi^{\tau(p)}\bigr)^*{\rm d}\lambda\bigr)_p(v_1,v_2) \\
&= {\rm d}\lambda_p(v_1,v_2) = \omega_p(v_1,v_2),
\end{aligned} \]
where it is used that $(\phi^t)^*\lambda=\lambda$ holds for all $t\in\R$, and that ${\rm d}\lambda(X,\cdot)\equiv0$. We obviously have $C \subset R$ and $C = \psi(C)$. After perhaps shrinking $S$, it can be assumed that $S$ is an annulus.\footnote{The case of a M\"obius band is excluded since $S$ is orientable.}

\begin{Lemma}
\label{lemma_symp_coord}
One finds $\epsilon>0$ and coordinates $(x,y) \in \R/\Z \times (-\epsilon,\epsilon)$ defined on some neighborhood of $C$ in $S$ such that $C = \R/\Z \times \{0\}$ and $\omega = {\rm d}x \wedge {\rm d}y$.
\end{Lemma}

\begin{proof}
Choose a smooth global periodic coordinate $x \in \R/\Z$ on $C$, and choose a smooth vector field~$Z$ which is transverse to $\Sigma$, tangent to~$S$, and such that ${\rm d}\lambda(\partial_x,Z) \equiv 1$ along~$C$. Denote its flow by $\Phi^t$. We find $\epsilon_0>0$ such that $\Phi^t(p)$ is defined for all $(p,t)\in \Sigma \times [-\epsilon_0,\epsilon_0]$. Then
\[
 S_0 = \bigl\{ \Phi^y(p) \mid (p,y) \in C \times (-\epsilon_0,\epsilon_0) \bigr\}
  \]
is an open neighborhood of~$C$ in~$S$ which is equipped with smooth coordinates $(x,y) \in \R/\Z \times (-\epsilon_0,\epsilon_0)$. In these coordinates $C = \R/\Z \times \{0\}$ and~$\omega$ agrees with the $2$-form ${\rm d}x \wedge {\rm d}y$ along $\R/\Z \times \{0\}$.

We now apply Moser's trick. Set $\omega_0 = {\rm d}x \wedge {\rm d}y$, $\omega_s = (1-s)\omega_0 + s\omega$, $s\in[0,1]$. In particular, every $\omega_s$ is an area-form coinciding with $\omega_0$ along $\R/\Z \times \{0\}$. Note that $\omega-\omega_0$ is closed and vanishes along ${\R/\Z \times \{0\}}$. We find a $1$-form $\alpha$ defined on $\R/\Z \times (-\epsilon_0,\epsilon_0)$ that vanishes on ${\R/\Z \times \{0\}}$ and satisfies ${\rm d}\alpha = \omega-\omega_0$. In fact, if $\omega-\omega_0 = f(x,y) \, {\rm d}x\wedge {\rm d}y$, then set ${\alpha = a(x,y) \, {\rm d}x}$, where $a(x,y) = -\int_0^y f(x,w){\rm d}w$. Consider the $s$-dependent vector field $W_s$ defined on $\R/\Z \times (-\epsilon_0,\epsilon_0)$ by the equation $i_{W_s}\omega_s + \alpha = 0$. Note that $W_s$ is well defined since $\omega_s$ is an area-form, and that $W_s$ vanishes on $\R/\Z \times \{0\}$ for all $s$. Hence, we find $\epsilon \leq \epsilon_0$ such that the initial value problem $\frac{{\rm d}}{{\rm d}s}\varphi^s=W_s \circ \varphi^s$, $\varphi^0 = {\rm id}$, defines a smooth family of embeddings
\[
\varphi^s\colon \ \R/\Z \times (-\epsilon,\epsilon) \to \R/\Z \times (-\epsilon_0,\epsilon_0)
\]
for all $s\in[0,1]$. One has $\varphi^s(x,0) = (x,0)$ for all $s \in [0,1]$ and $x\in\R/\Z$. Differentiate in $s$ to obtain
\[
\frac{{\rm d}}{{\rm d}s} (\varphi^s)^*\omega_s = (\varphi^s)^*({\rm d}i_{W_s}\omega_s + {\rm d}\alpha) = {\rm d}(\varphi^s)^*(i_{W_s}\omega_s + \alpha) = 0.
\]
In particular, \smash{$\bigl(\varphi^1\bigr)^*\omega = \omega_0$} and one obtains the desired coordinates by replacing $x$ and $y$ by \smash{$x\circ \bigl(\varphi^1\bigr)^{-1}$} and \smash{$y\circ \bigl(\varphi^1\bigr)^{-1}$}, respectively.
\end{proof}

Let $(x,y) \in \R/\Z \times (-\epsilon,\epsilon)$ be the $\omega$-symplectic coordinates given by Lemma~\ref{lemma_symp_coord}, defined on a neighborhood of~$C$ in~$S$. We can assume that this neighborhood is contained in the neighborhood~$R$ of~$C$ in~$S$. Hence, we get from~\eqref{return_map} a well-defined map
\begin{equation}
\label{local_return_map_in_coordinates}
\psi\colon \ \R/\Z\times(-\epsilon,\epsilon) \to \R/\Z \times \R
\end{equation}
that preserves ${\rm d}x \wedge {\rm d}y$, leaves $C = \R/\Z \times \{0\}$ invariant, and defines an embedding onto an open subset of the target. The diffeomorphism~$\psi|_C$ is necessarily orientation-preserving: otherwise it would have fixed points, in contradiction to~(P3). Hence, if we write the return map $\psi$ in components
\begin{equation}
\label{components_local_return_map}
\psi(x,y) = (X(x,y),Y(x,y)),
\end{equation}
we get
\begin{equation}
\label{signs_partial_der}
D_1X(x,0)>0, \qquad D_2Y(x,0)>0 \qquad \forall x\in\R/\Z.
\end{equation}

\begin{Lemma}
\label{lemma_extension}
There exists $r \in (0,\epsilon)$ and a smooth map
\[
 \varphi \colon \ \R/\Z \times (-r,+\infty) \to \R/\Z \times \R
 \] satisfying:
\begin{itemize}\itemsep=0pt
\item $\varphi$ coincides with $\psi$ on $\R/\Z \times (-r,r)$.
\item $\varphi$ preserves ${\rm d}x \wedge {\rm d}y$.
\item $\varphi$ maps $\R/\Z\times[0,+\infty)$ diffeomorphically onto itself.
\end{itemize}
\end{Lemma}

\begin{proof}
The proof is performed by considering a generating function representing the map $\psi$ on a small neighborhood of $\R/\Z\times \{0\}$, and then carefully extending it to an open neighborhood of $\R/\Z\times[0,+\infty)$ in order to define the map~$\varphi$.

In view of~\eqref{signs_partial_der}, the map $h\colon (x,y) \mapsto (x,Y)$ defines a diffeomorphism between some open neighborhood of $C=\R/\Z\times\{0\}$ and its image. It fixes every point of~$C$. Choose $\delta_1>0$ such that $\R/\Z\times (-\delta_1,\delta_1)$ is contained in the image of~$h$. Consider the open neighborhood $U_1 = h^{-1}(\R/\Z\times (-\delta_1,\delta_1))$ of $\R/\Z \times\{0\}$. It then follows that $h|_{U_1}\colon U_1\to \R/\Z\times (-\delta_1,\delta_1)$ is a~diffeomorphism. On $\R/\Z\times (-\delta_1,\delta_1)$, we can consider $(x,Y)$ as independent variables and write~$X$, $y$ as functions of $(x,Y)$ defined on this set. On $\R/\Z\times (-\delta_1,\delta_1)$ we have well-defined maps $f(x,Y) = (x,y)$, $g(x,Y) = (X,Y)$. Note that $f$ is the inverse of $h|_{U_1}$. Since ${\rm d}X \wedge {\rm d}Y = {\rm d}x \wedge {\rm d}y$, $\eta = (y-Y){\rm d}x + (X-x){\rm d}Y$ is a closed $1$-form on $\R/\Z\times(-\delta_1,\delta_1)$. Its period vanishes since $y(x,0)=0$ holds for all~$x$. Here we used that $(x,0)\mapsto \psi(x,0)$ is an orientation preserving diffeomorphism of~$C$. We find a smooth function $S=S(x,Y)$ on $\R/\Z\times (-\delta_1,\delta_1)$ such that
\begin{equation*}
X = x + D_2S(x,Y), \qquad y = Y + D_1S(x,Y).
\end{equation*}
This means $S$ is a generating function for $\psi$. We have that $D_1S(x,0)=0$ and $D_{12}S(x,0)>-1$ hold for all $x\in\R/\Z$. Up to shrinking~$\delta_1$, it can further be assumed without loss of generality that
\begin{equation}\label{estimate_mixed_at_C}
\alpha = \inf_{\R/\Z\times(-\delta_1,\delta_1)}  D_{12}S   >  -1
\end{equation}
and that $\psi|_{U_1} = g \circ h = g \circ f^{-1}$. We proceed assuming that $\delta_1$ meets these requirements. The set $U_1$ is shrunk accordingly. We can now fix $r>0$ from the statement: choose it so that $\R/\Z\times (-r,0] \subset U_1$.

Fix a smooth function $\beta\colon \R \to [0,1]$ satisfying $\beta'(s) \leq 0$, $\beta(s) = 1$ for all $s$ on a neighborhood of $(-\infty,0]$, and ${\rm supp}(\beta) \subset (-\infty,1)$. For every $\delta \in (0,\delta_1)$, consider the smooth function $\beta_\delta(s) = \beta(s/\delta)$. In order to avoid ambiguity of notation, consider independent variables $(x,v) \in \R/\Z\times (-\delta_1,+\infty)$. Let $S_\delta\colon \R/\Z \times (-\delta_1,+\infty) \to \R$ be given~by
\begin{equation*}
S_\delta(x,v) =   \begin{cases}  \beta_\delta(v)S(x,v) + (1-\beta_\delta(v))D_2S(x,0)v & \text{if $v \in (-\delta_1,\delta_1)$}, \\D_2S(x,0)v  & \text{if $v\geq \delta_1$}.
 \end{cases}
\end{equation*}
Note that the functions $S_\delta$ and $S$ coincide on a neighborhood of $\R/\Z \times (-\delta_1,0]$ in $\R/\Z \times (-\delta_1,+\infty)$.

We claim that if $\delta$ is small enough, then $S_\delta$ is a generating function for a map meeting the requirements of the statement. We argue as follows. Define functions $X_\delta(x,v)$, $y_\delta(x,v)$ by
\begin{equation*}
X_\delta(x,v) = x + D_2S_\delta(x,v), \qquad y_\delta(x,v) = v + D_1S_\delta(x,v).
\end{equation*}
Consider smooth maps defined by
\begin{equation}
\label{partial_maps}
f_\delta\colon \ (x,v) \mapsto (x,y_\delta(x,v)), \qquad g_\delta\colon \ (x,v) \mapsto (X_\delta(x,v),v).
\end{equation}
All functions $X_\delta$, $y_\delta$, $f_\delta$, $g_\delta$ are defined on $\R/\Z \times (-\delta_1,+\infty)$.
The derivative of $f_\delta$ is
\begin{equation*}
\begin{bmatrix}
1 & 0 \\
D_{11}S_\delta(x,v) & 1 + D_{12}S_\delta(x,v)
\end{bmatrix}
\end{equation*}
while the derivative of $g_\delta$ is
\begin{equation*}
\begin{bmatrix}
1 + D_{21}S_\delta(x,v) & D_{22}S_\delta(x,v) \\
0 & 1
\end{bmatrix}.
\end{equation*}
Both matrices above are invertible provided $D_{21}S_\delta(x,v) = D_{12}S_\delta(x,v) > -1$. Let us now check this inequality. On $\R/\Z \times (-\delta_1,\delta_1)$, simple calculations using $D_1S(x,0)\equiv0$ lead to
\begin{equation}\label{expression_mixed_partial}
D_{12}S_\delta(x,v) = D_{12}S(x,0) + \beta_\delta'(v)v \Delta_1(x,v) + \beta_\delta(v) \Delta_2(x,v),
\end{equation}
where
\begin{align*}
& \Delta_1(x,v)  = \frac{D_1S(x,v)-D_1S(x,0)}{v} - D_{21}S(x,0), \qquad
 \Delta_2(x,v)= D_{12}S(x,v) - D_{12}S(x,0).
\end{align*}
Note that \[
 \lim_{v\to0}   \sup_{x\in\R/\Z}  ( |\Delta_1(x,v)| + |\Delta_2(x,v)|  ) = 0.
\]
Moreover, $\beta_\delta'(v)v$, $\beta_\delta(v)$ vanish for $v \geq \delta$ and are bounded from above by a constant independent of $\delta$. It then follows from these facts and~\eqref{expression_mixed_partial} that
\[ \lim_{\delta\to0^+}  \sup_{(x,v) \in \R/\Z \times [0,\delta]}   |D_{12}S_\delta(x,v)-D_{12}S(x,0)| = 0. \]
In particular, we get from~\eqref{estimate_mixed_at_C} and the above formula that there exists $\alpha_1 \in (-1,\alpha]$ such that $D_{21}S_\delta(x,v) \geq \alpha_1$ holds for all $(x,v) \in \R/\Z \times [0,\delta]$ provided $\delta$ is small enough. If $v > \delta$, then $S_\delta(x,v) = D_2S(x,0)v$ and $D_{21}S_\delta(x,v) \geq \alpha > -1$ follows from~\eqref{estimate_mixed_at_C}. If $v \in (-\delta_1,0]$, then $D_{21}S_\delta(x,v)$ coincides with $D_{21}S(x,v) \geq \alpha > -1$. We have shown that there exists $\delta_* \in (0,\delta_1)$ such that
\begin{equation}
\label{est_mixed_with_delta}
D_{21}S_\delta(x,v) = D_{12}S_\delta(x,v) \geq \alpha_1 > -1
\end{equation}
holds on $\R/\Z \times (-\delta_1,+\infty)$ provided $\delta \in (0,\delta_*)$. In particular, both maps in~\eqref{partial_maps} have non-singular derivatives on $\R/\Z \times (-\delta_1,+\infty)$.

We now argue to show that both maps in~\eqref{partial_maps} map $\R/\Z\times[0,+\infty)$ diffeomorphically onto itself provided $\delta \in (0,\delta_*)$. We claim that for each fixed $x\in\R/\Z$ and $\delta \in (0,\delta_*)$ the map $v\mapsto y_\delta(x,v)$ is a self-diffeomorphism of $[0,+\infty)$. This follows from $y_\delta(x,0)=0$ and~\eqref{est_mixed_with_delta} since $\frac{{\rm d}}{{\rm d}v}y_\delta(x,v) = 1+D_{21}S_\delta(x,v)$. As a~consequence, one obtains that $f_\delta$ defines a bijection of $\R/\Z \times [0,+\infty)$ onto itself: given $(\bar x,\bar y) \in \R/\Z \times [0,+\infty)$, the equation $f_\delta(x,v) = (\bar x,\bar y)$ has a~unique solution, where $x = \bar x$ and $v$ is uniquely determined by $y_\delta(\bar x,v)=\bar y$. Note that $\det Df_\delta > 0$ everywhere. Hence, $f_\delta$ maps $\R/\Z\times[0,+\infty)$ diffeomorphically onto itself provided $\delta \in (0,\delta_*)$. The argument for $g_\delta$ is similar: \eqref{est_mixed_with_delta} implies that $x \mapsto X_{\delta}(x,v)$ is an orientation-perserving self-diffeomorphism of $\R/\Z$ for all $v\geq 0$.

Since $S$ and $S_\delta$ coincide on a neighborhood of $\R/\Z\times(-\delta_1,0]$ in $\R/\Z\times(-\delta_1,+\infty)$, the maps $f_\delta$, $g_\delta$ extend the maps $f,g$ from such a neighborhood to $\R/\Z\times(-\delta_1,+\infty)$, respectively. Moreover, by what is shown above, they define diffeomorphisms onto their images. Note that $\R/\Z \times (-r,+\infty)$ is contained in the image of $f_\delta$. The proof is complete if we set $\varphi = g_\delta \circ (f_\delta)^{-1}$ with any arbitrarily chosen $\delta \in (0,\delta_*)$. The map $\varphi$ preserves the area form since it has $S_\delta$ as a~generating function.
\end{proof}

Consider the diffeomorphism
\begin{equation*}
h\colon \ \C\setminus\{0\} \to \R/\Z \times (0,+\infty), \qquad h\bigl(u+{\rm i}v=r{\rm e}^{{\rm i}\theta}\bigr) = \bigl(-\theta/2\pi,\pi r^2\bigr),
\end{equation*}
which is area-preserving
\begin{equation*}
h^*({\rm d}x \wedge {\rm d}y) = {\rm d}u \wedge {\rm d}v
\end{equation*}
as one can check. Consider the diffeomorphism $f\colon \C\setminus\{0\}\to \C\setminus\{0\}$ defined by $f = h^{-1} \circ \varphi \circ h$, where $\varphi$ is given by Lemma~\ref{lemma_extension}. Then~$f$ preserves orientations and Lebesgue measure, and satisfies $f(z) \to 0$ as $z\to0$. Thus, it can be extended to an orientation-preserving homeomorphism $f\colon \C\to\C$ by requiring $f(0)=0$.

\begin{Lemma}\label{lemma_0_is_isolated}
There exists a neighborhood $V$ of $0$ in $\C$ such that $\{0\} = \bigcap_{k\in\Z}f^k(V)$.
\end{Lemma}

\begin{proof}
Recall that $h \circ f \circ h^{-1}$ agrees with the map $\psi$ from~\eqref{local_return_map_in_coordinates} on $\R/\Z\times (0,b)$ provided $b>0$ is small enough. Hence if the lemma is not true, then we find for each $n \in \N$ a number $a_n>0$ and a point $(x_n,y_n) \in \R/\Z\times(0,+\infty)$ such that $a_n \to 0$, $\psi^k(x_n,y_n)$ is well defined and $\psi^k(x_n,y_n) \in \R/\Z \times (0,a_n)$, for every $(n,k)$ in $\N\times\Z$. The points $p_{n,k} := \psi^k(x_n,y_n)$ correspond to points in $S$ uniformly and arbitrarily close to the invariant torus $\Sigma$. Moreover, for each $n$ the set $\{p_{n,k}\}_{k\in\Z}$ is contained in a single trajectory of the Reeb flow~$\phi^t$: In fact, $p_{n,k+1} = \phi^{\tau(p_{n,k})}(p_{n,k})$ holds for all $(n,k)$. In particular, $\{p_{n,k}\}_{k\in\Z} \subset \phi^\R(p_{n,0})$. Here $\tau$ is the return-time function~\eqref{return_time_function}. The continuity of the flow together with the fact that~$\tau$ is bounded on a neighborhood of~$C$ in~$S$ directly imply that for every neighborhood $W$ of $\Sigma$ there exists $n_0$ such that
\[
 n\geq n_0 \Rightarrow \phi^{[0,\tau(p_{n,k})]}(p_{n,k}) \subset W \setminus \Sigma,
 \]
in particular
\[ n\geq n_0 \Rightarrow \phi^\R(p_{n,0}) \subset W \setminus \Sigma.
\] This is in contradiction to assumption (P2) asserting that $\Sigma$ is isolated with respect to the Reeb flow.
\end{proof}

By the previous lemma, the Lefschetz index ${\rm ind}(f^j,0)$ of $0$ as a fixed point of $f^j$ is well defined, for every $j\in\N$.

\begin{Lemma}\label{lemma_index_is_positive}
The identity ${\rm ind}(f^j,0)=1$ holds for each $j\in\N$.
\end{Lemma}

\begin{proof}
Consider the map $\psi$ in components $\psi(x,y) = (X(x,y),Y(x,y))$ as in~\eqref{components_local_return_map}. Then on a~punctured neighborhood of~$0$, we write
\[
f\bigl(r{\rm e}^{{\rm i}\theta}\bigr) = R{\rm e}^{{\rm i}\vartheta}, \qquad \text{where} \quad (r,\theta) = \bigl(\sqrt{y/\pi},-2\pi x\bigr), \quad (R,\vartheta) = \bigl(\sqrt{Y/\pi},-2\pi X\bigr).
\]
For each $\lambda\in(0,1]$, consider the map
\[
f_\lambda(z) = \frac{1}{\lambda}f(\lambda z).
\]
Then $f=f_1$ and each $f_\lambda$ is an area-preserving homeomorphism of $\C$ which restricts to a smooth self-diffeomorphism of $\C\setminus\{0\}$. Note that $(f^j)_\lambda = (f_\lambda)^j$.

The first task is to understand the behavior of $f_\lambda$ as $\lambda\to0^+$. From the above formulas, we~get
\begin{align*}	
\lim_{\lambda \to 0^+} f_\lambda\bigl(r{\rm e}^{{\rm i}\theta}\bigr) &= \lim_{\lambda\to0^+}\frac{1}{\lambda} \sqrt{\frac{Y(x,\lambda^2y)}{\pi}} {\rm e}^{-{\rm i}2\pi X (x,\lambda^2y)} \\
&  = \lim_{\lambda\to0^+}\sqrt{\frac{D_2Y(x,0)y}{\pi} + \frac{o(\lambda^2|y|)}{\lambda^2}} {\rm e}^{-{\rm i}2\pi X(x,\lambda^2y)}
 = rg(\theta){\rm e}^{{\rm i}(\theta+\Delta(\theta))} =: f_0\bigl(r{\rm e}^{{\rm i}\theta}\bigr),
\end{align*}
where
\begin{equation*}
g(\theta) = \sqrt{D_2Y(-\theta/2\pi,0)}, \qquad \Delta(\theta) = -2\pi(X(-\theta/2\pi,0)+\theta/2\pi).
\end{equation*}
The convergence is in $C^\infty_{\rm loc}(\C\setminus \{0\})$ and in $C^0_{\rm loc}(\C)$. Note that $g(\theta)>0$ for all $\theta$.

Since the rotation number of the map $x \in \R/\Z \mapsto X(x,0) \in \R/\Z$ is irrational, we can be sure that $\Delta(\theta) \not\in 2\pi\Z$, for all $\theta$. In particular, $f_0$ has no fixed points on $\C\setminus\{0\}$ and ${\rm ind}(f_0,0)$ is well defined. By continuity of the degree under continuous homotopies, we know that
\begin{equation*}
{\rm ind}(f,0) = {\rm ind}(f_0,0).
\end{equation*}
We will now show that ${\rm ind}(f_0,0)=1$. This Lefschetz index is the winding number of $f_0\bigl(r{\rm e}^{{\rm i}\theta}\bigr)-r{\rm e}^{{\rm i}\theta}$ around zero, where $r$ is small and $\theta$ increases from~$0$ to~$2\pi$. Write
\begin{equation*}
f_0\bigl(r{\rm e}^{{\rm i}\theta}\bigr)-r{\rm e}^{{\rm i}\theta} = r{\rm e}^{{\rm i}\theta}\bigl( g(\theta){\rm e}^{{\rm i}\Delta(\theta)} - 1 \bigr)
\end{equation*}
to conclude that ${\rm ind}(f_0,0) = m+1$, where $m\in\Z$ is the winding number of $g(\theta){\rm e}^{{\rm i}\Delta(\theta)}$ around $1\in\C$. The preservation of area tells us that
\begin{equation*}
\begin{aligned}
r {\rm d}r\wedge {\rm d}\theta &= (rg(\theta)) {\rm d}(rg(\theta)) \wedge {\rm d}(\theta+\Delta(\theta)) \\
&= (rg(\theta)) (g(\theta){\rm d}r + rg'(\theta){\rm d}\theta) \wedge (1+\Delta'(\theta)){\rm d}\theta
 = g(\theta)^2(1+\Delta'(\theta))   r{\rm d}r\wedge {\rm d}\theta
\end{aligned}
\end{equation*}
which is then equivalent to
$g(\theta)^2(1+\Delta'(\theta)) = 1$.
In particular,
\begin{equation*}
\begin{aligned}
& g(\theta)<1 \Leftrightarrow 1+\Delta'(\theta)>1 \Leftrightarrow \Delta'(\theta)>0, \\
& g(\theta)>1 \Leftrightarrow 1+\Delta'(\theta)<1 \Leftrightarrow \Delta'(\theta)<0.
\end{aligned}
\end{equation*}
Geometrically it is clear from the above that $m\leq 0$: $g(\theta){\rm e}^{{\rm i}\Delta(\theta)}$ can only wind around~$1$ in the clockwise direction. This is so because $m$ is equal to the intersection number of this curve with the ray $(1,+\infty)$ equipped with the orientation from $1$ to $+\infty$. Since $g(\theta)>1 \Rightarrow \Delta'(\theta)<0$, every intersection counts negatively. If $m<0$, then $g(\theta){\rm e}^{{\rm i}\Delta(\theta)}$ hits $(1,+\infty)$, which then implies that $\Delta(\theta)$ hits $2\pi\Z$, a contradiction to irrationality of the rotation number. We conclude that
\begin{equation*}
m=0 \qquad \text{and} \qquad {\rm ind}(f_0,0) = 1
\end{equation*}
as claimed before.

Now note that for every $j\geq1$ the map $f_0^j$ has the same form as $f_0$:
\[
 f_0^j\bigl(r{\rm e}^{{\rm i}\theta}\bigr) = rg_j(\theta){\rm e}^{{\rm i}(\theta+\Delta_j(\theta))}
\]
for functions $g_j(\theta)$, $\Delta_j(\theta)$ with the same properties as before. In particular, $\Delta_j(\theta)$ never hits $2\pi\Z$, by irrationality of the rotation number. The same arguments as above will show that ${\rm ind}\bigl(f_0^j,0\bigr)=1$.
\end{proof}

Now we can use Theorem~\ref{thm_LeCY} to achieve a contradiction and finish the proof of Theorem~\ref{thm_main_1}. The argument is as follows. Lemma~\ref{lemma_0_is_isolated} states that $0$ is a fixed point of~$f$ which is isolated as an invariant set of~$f$. Lemma~\ref{lemma_index_is_positive} says that ${\rm ind}(f^j,0)$ is equal to $1$ for all positive integers $j$, but Theorem~\ref{thm_LeCY} states that there are infinitely many positive integers $j$ such that ${\rm ind}(f^j,0)$ is non-positive. This is absurd and we conclude that the hypothesis that $\Sigma$ contains no periodic orbits cannot hold.

\subsection*{Acknowledgements}

We are grateful to the referees for the careful read of our arguments.
UH is partially funded by the DFG-Graduiertenkolleg \textit{Energy, Entropy, and Dissipative Dynamics $($EDDy$)$}, by the FAPERJ grant E02/2025 \textit{Pesquisador visitante}, and by the DFG-funded SFB/TRR 191 \textit{Symplectic Structures in Geometry, Algebra and Dynamics}.

\pdfbookmark[1]{References}{ref}
\LastPageEnding

\end{document}